\documentclass[reqno,12pt]{amsart}

\RequirePackage{amssymb}
\RequirePackage{times}
\RequirePackage{graphicx}
\RequirePackage{hyperref}
\RequirePackage{lipsum}
\RequirePackage[normalem]{ulem}
\RequirePackage[hang,small,bf]{caption}

\newcommand{\R}{\mathbb{R}}

\newcommand{\E}{\mathbb{E}}
\newcommand{\bP}{\mathbb{P}}

\newcommand{\df}{\stackrel{{{\mathrm{def}}}}{=}} 

\newtheorem{theorem}{Theorem}[section]
\newtheorem{lemma}[theorem]{Lemma}

\newtheorem{corollary}[theorem]{Corollary}

\makeatletter
\g@addto@macro{\endabstract}{\@setabstract}
\newcommand{\authorfootnotes}{\renewcommand\thefootnote{\@fnsymbol\c@footnote}}%
\makeatother

\makeatletter
\let\@fnsymbol\@arabic
\makeatother

\begin{document}
\begin{center}
  \LARGE 
  On estimation in the reduced-rank 
regression with a large number of responses and predictors \par \bigskip
  \normalsize
  \authorfootnotes
  Vladislav Kargin\footnote{Email: vladislav.kargin@gmail.com; current address: 282 Mosher Way, Palo Alto, CA 94304, USA} \par \bigskip

	\bigskip

\end{center}

\begin{center}
\textbf{Abstract}
\end{center}

\begin{quotation}
We consider a multivariate linear response regression in which the number of
responses and predictors is large and comparable with the number of
observations, and the rank of the matrix of regression coefficients is assumed to be
small. We study the distribution of singular values for the matrix of regression
coefficients and for the matrix of predicted responses. For both matrices, it is found that the
limit distribution of the largest singular value is a rescaling of the Tracy-Widom
distribution. Based on this result, we suggest algorithms for the model rank
selection and compare them with the algorithm suggested by Bunea, She and
Wegkamp. Next, we design two consistent estimators for the singular values of
the coefficient matrix, compare them, and derive the asymptotic distribution for one of these 
estimators.
\end{quotation}

\section{Introduction}

In this paper we are concerned with the reduced rank variant of the
multivariate response regression model. We are given $N$ observations of the
predictors $X_{i}\in \mathbb{R}^{p}$ and responses $Y_{i}\in \mathbb{R}^{r}, 
$ which are assumed to be related by the linear regression model: 
\begin{equation}
Y=XA+U,  \label{regression_model}
\end{equation}%
where $A$ is an unknown $p$-by-$r$ matrix and $U$ is a noise matrix. This model is ubiquitous in statistics, signal
processing, and numerical analysis.

On methodological grounds one often postulates that
the responses depend only on a small number of factors which are linear
combinations of the predictors. This postulate leads to a 
model, in which $A$ is assumed to be a low-rank matrix: 
\begin{equation}
A=\sum_{j=1}^{s}\theta _{j}u_{j}v_{j}^{\ast },  \label{low_rank_model}
\end{equation}%
where $\left\{ u_{j}\in \mathbb{R}^{p}\right\} $ and $\left\{ v_{i}\in 
\mathbb{R}^{r}\right\} $ are two fixed orthonormal vector systems. This
model appeared already in Anderson (1951) \cite{anderson51}, and it was named \textit{reduced-rank regression%
} in Izenman
(1975) \cite{izenman75}.  In some contexts, this model is also known under the names
\textit{simultaneous linear prediction}  (Fortier (1966) 
\cite{fortier66}) and \textit{redundancy analysis} (van den Wollenberg
(1977) \cite{wollenberg77}), both of which assume that $U$ has the covariance matrix equal to $\sigma^2 I$.
The reduced-rank model has been intensively studied, and many results are collected in the
monograph by Reinsel and Velu (1998) \cite{reinsel_velu98}.

In this paper, we assume that $U$ has the covariance matrix equal to $\sigma^2 I$, and we are interested in the situation in which all three
variables, $p,$ $r,$ and $N,$ grow at the same rate.

\textbf{Assumption A1.} \textit{It is assumed that  as $N\rightarrow
\infty$,  $\frac{N}{p}\rightarrow $ $%
1+\lambda \geq 1$ and $\frac{N}{r}\rightarrow \mu >0$.}

It is also useful to define $\beta \df \lim_{N\rightarrow \infty }\frac{p}{r}=\mu
/\left( 1+\lambda \right).$

The studies devoted to
the reduced-rank regression in this setup are relatively recent and include
Bunea, She, and Wegkamp (2011) \cite{bunea_she_wegcamp11} and Giraud (2011) 
\cite{giraud11}.

We address the following questions. First, is it possible to detect that the
true matrix $A$ is not zero? If yes, then how do we how do we estimate the rank and singular values of $A$?

Our approach to these questions is based on the study of the statistical
properties of the standard least squares estimator 
\begin{equation*}
\widehat{A}\df X\backslash Y\equiv \left( X^{\ast }X\right) ^{-1}X^{\ast }Y
\end{equation*}%
and the matrix of fitted responses: 
\begin{equation*}
\widehat{Y}\df X\widehat{A}
\end{equation*}

By using this approach, we will develop a rank-selection algorithm 
which
performs better than the algorithm from \cite{bunea_she_wegcamp11} in a
certain range of parameters and is simpler than the algorithm in \cite{yelm07}. In
addition, we will develop tools for consistent estimation of singular values $\theta
_{i}.$ The paper \cite{bunea_she_wegcamp11} does not address this issue,
since its focus is on minimizing the prediction error, in particular on
bounds for $\E\left\Vert XA-X\widetilde{A}\right\Vert ,$ where $\widetilde{A}$
is an estimator of $A$ and the expectation is over randomness in $U$.

The rest of the paper is organized as follows. Section \ref{section_major_results} describes the major results. Section \ref%
{subsection_largest_eigenvalue_nonnull} provides the details of the proofs. Section \ref{conclusion} recapitulates the results.  And Appendix \ref%
{section_limiting_distribution} provides a proof for the theorem about the limiting
distribution of singular values of $\widehat{A}$.

\section{Major results}
\label{section_major_results}

\subsection{Tests of the null hypothesis}

Let $X$ be a $p$-by-$r$ real Gaussian matrix: each row is an independent observation from $\mathcal{N}(0,\Sigma)$. Then, an $r$-by-$r$ matrix $X^{\ast}X$ is said to be a \emph{Wishart matrix} with distribution $W_r(\Sigma,p)$.

A random $m$-by-$m$ matrix $X$ is said to belong to the \emph{(real) Jacobi ensemble} with parameters $\alpha_1$ and $\alpha_2$, if its distribution is invariant with respect to orthogonal transformations and the distribution of its eigenvalues is given by 
  
	\begin{equation}
f^{(\alpha_1,\alpha_2)}\left( \lambda_{1},\ldots ,\lambda_{m}\right) =\frac{1}{c}\prod\limits_{j=1}^{m}
\lambda_{j}^{\alpha_1}( 1-\lambda_i) ^{\alpha_2}
\prod\limits_{1\leq j<k\leq m}\left\vert
\lambda_{j}-\lambda_{k}\right\vert.  \label{pdf_Jacobi}
\end{equation}%

The following result is fundamental for the analysis of matrices $\widehat{A}
$ and $\widehat{Y}.$

  \begin{theorem}
	\label{theorem_fundamental}
	(i) Suppose that $U$ is an $N$-by-$r$ matrix with i.i.d standard real Gaussian entries, and $X$ is an  $N$-by-$p$ full-rank matrix ($N\geq p$) independent of $U$. Then the squared singular values of $\widehat{Y}\df X(X\backslash U)$ are distributed as the eigenvalues of the Wishart matrix with distribution $W_r(I,p)$. \\ 
	(ii) In addition, suppose that $X$ has i.i.d standard real Gaussian entries.
	Let $s^2_i$ be the squared singular values of $\widehat{A}\df X\backslash U$ and $f_i=s^2_i/(1+s^2_i)$. Then, the positive $f_i$ are distributed as eigenvalues of 
			the Jacobi ensemble with parameters $m=\min\{p,r\}$, $\alpha_1=(|r-p|-1)/2$ and $\alpha_2=(N-p-1)/2$.  
  \end{theorem} 
	\textbf{Proof}: 	The matrix $ \widehat{Y}=X(X\backslash U)$ is the orthogonal projection of $r$ column vectors of $U$ on the 
	$p$-dimensional column span of $X$. Hence, in an
	appropriate basis, $\widehat{Y}$ is a block matrix with one block given by a $p$-by-$r$ matrix with i.i.d. standard Gaussian entries and another block of $N-p$-by-$r$ matrix of zeros. This proves the first part of the theorem. 	For the second part, note that positive eigenvalues of $\widehat{A}^{\ast} \widehat{A} = U^{\ast }X\left( X^{\ast }X\right) ^{-2}X^{\ast }U$
have the same distribution as positive eigenvalues of $B^{-1}C,$ where $B$ and $C$ are independent Wishart matrices. 

Indeed, the rank of matrices  $U^{\ast}X(X^{\ast}X)^{-2}X^{\ast}U$ and \\ $X(X^{\ast}X)^{-2}X^{\ast}UU^{\ast}$ is 
$\min\{p,r\}$, and their positive eigenvalues are the same. Let $W$ be an orthogonal $N$-by-$p$ matrix formed by the eigenvectors of $X(X^{\ast}X)^{-2}X^{\ast}$ and such that the matrix $W^{\ast}X(X^{\ast}X)^{-2}X^{\ast}W$ is diagonal with positive eigenvalues on the diagonal.
  These eigenvalues coincide with positive eigenvalues of the inverse of a Wishart matrix, $(X^{\ast}X)^{-1}$, where the Wishart matrix has the distribution $W_p(I,N)$. The matrix $W^{\ast}UU^{\ast}W$ is Wishart with distribution $W_p(I,r)$. 
 
In addition, matrices  $W^{\ast}X(X^{\ast}X)^{-2}X^{\ast}W$ and $W^{\ast}UU^{\ast}W$ are independent because the eigenvalues and eigenvectors of $X(X^{\ast}X)^{-2}X^{\ast}$ are independent.  Finally, since similarity transformations do not change eigenvalues, the distribution of positive eigenvalues of  $X(X^{\ast}X)^{-2}X^{\ast}UU^{\ast}$ is the same as the distribution of positive eigenvalues of $B^{-1}C$ for two independent Wishart matrices $B$ and $C$ with distributions $W_p(I,N)$ and $W_p(I,r)$, respectively.

Next, note that the eigenvalues of $B^{-1}C$ are the same as those
of $F\left( I-F\right) ^{-1},$ where $F=(B+C)^{-1}C.$ Hence the
eigenvalues of $B^{-1}C,$ denoted by $l_{i},$ are
related to the eigenvalues of $\left( B+C\right) ^{-1}C,$ denoted by $f_i$ by the transformation,%
\begin{equation}
l_{i}=\frac{f_{i}}{1-f_{i}}.  \label{relation_of_eigenvalues}
\end{equation}%
Then one can use the classical fact that the eigenvalues of $\left( B+C\right) ^{-1}C$
are distributed as in (\ref{pdf_Jacobi}) with $\alpha_1=(|r-p|-1)/2$ and $\alpha_2=(N-p-1)/2$. (See Theorem 3.3.1 in Muirhead 
\cite{muirhead82}).
\hfill $\square$

There are several ways to test the null hypothesis that $A=0.$ The simplest
way is to compute the largest eigenvalue of $Y^{\ast }Y,$ called $l_{Y,1},$
and compute 
\begin{equation*}
s_{Y,1}=\frac{l_{Y,1}-\mu ^{\left( 1\right) }}{\sigma ^{\left(
1\right) }},
\end{equation*}%
where 
\begin{eqnarray}
\mu ^{\left( 1\right) } &=&\left( \sqrt{N-1}+\sqrt{r}\right) ^{2},
\label{test1_scaling} 
\\
\sigma ^{\left( 1\right) } &=&\left( \sqrt{N-1}+\sqrt{r}\right) \left[ \frac{%
1}{\sqrt{N-1}}+\frac{1}{\sqrt{r}}\right] ^{1/3}.  \notag
\end{eqnarray}

If $A=0$ and the noise matrix $U$ has i.i.d. standard Gaussian entries,  then by
Theorem 1.1 in Johnstone (2001) \cite{johnstone01} $s_{Y,1}$ converges in distribution to the Tracy-Widom
distribution $\mathcal{F}_{1}.$\footnote{Some quantiles of the Tracy-Widom distribution are: $x_{50\%}=-1.3,$ $%
x_{10\%}=0.45,$ $x_{5\%}=0.98,$ $x_{2\%}=1.60,$ $x_{1\%}=2.02,$ where $%
\mathbb{P}\left( x\geq x_{c}\right) =c.$} Note that this test does not depend on  $X$ and that the normality assumption on the entries of $U$ can be significantly relaxed. Indeed, by results of Pillai and Yin \cite{pillai_yin14}, the convergence to the Tracy-Widom distribution holds provided that the entries are independent, have zero mean, unit variance, and subexponential decay: $\bP\{U_{ij}>t\}\leq \kappa^{-1}\exp\{-t^\kappa\}$ for a positive $\kappa$ and all $t>1$. (The first universality result of this type is due to Soshnikov \cite{soshnikov02}.)

Two other methods to test the null hypothesis are conditional on $X$. They are based on the singular values of
matrices $\widehat{A}\df X\backslash Y$ and $\widehat{Y}\df X \widehat{A}.$ Let us define two statistics, 
$s_{\widehat{Y},1}$ and $s_{\widehat{A},1}:$ 

\begin{equation*}
s_{\widehat{Y},1}\df\frac{l_{\widehat{Y},1}-\mu ^{\left( 2\right) }}{\sigma ^{\left(
2\right) }},
\end{equation*}%
where  $l_{\widehat{Y},1},$ is the square of the largest singular value of  $\widehat{Y}$ and 
\begin{eqnarray}
\mu ^{\left( 2\right) } &=&\left( \sqrt{p-1}+\sqrt{r}\right) ^{2},
\label{test2_scaling} 
\\
\sigma ^{\left( 2\right) } &=&\left( \sqrt{p-1}+\sqrt{r}\right) \left[ \frac{%
1}{\sqrt{p-1}}+\frac{1}{\sqrt{r}}\right] ^{1/3}.  \notag
\end{eqnarray}%

Similarly,
\begin{equation*}
s_{\widehat{A},1}=\frac{\log l_{\widehat{A},1}-\mu ^{\left( 3\right) }}{%
\sigma ^{\left( 3\right) }},
\end{equation*}%
where $l_{\widehat{A},1}$ is the square of the largest singular value of $\widehat{A}$ and 
\begin{eqnarray}
\mu ^{\left( 3\right) } &=&2\log \tan \left( \frac{\phi +\gamma }{2}\right) ,
\label{centering_and_scaling} 
\\
\sigma ^{\left( 3\right) } &=&\left[ \frac{16}{\left( N+r-1\right) ^{2}}%
\frac{1}{\sin ^{2}\left( \phi +\gamma \right) \sin \phi \sin \gamma }\right]
^{1/3},  \notag
\end{eqnarray}%
and the angle parameters $\gamma $ and $\phi $ are defined by 
\begin{eqnarray*}
\sin ^{2}\left( \frac{\gamma }{2}\right) &=&\frac{\min \left( p,r\right) -1/2%
}{N+r-1}, \\
\sin ^{2}\left( \frac{\phi }{2}\right) &=&\frac{\max \left( p,r\right) -1/2}{%
N+r-1}.
\end{eqnarray*}

\bigskip

\begin{theorem}
(i) Suppose that $A=0$, Assumption A1 is satisfied, $U$ is a matrix with
independent standard Gaussian entries, and $X$ is a full-rank matrix independent of $U.$
Then the random variable $s_{\widehat{Y},1}$ converges in distribution to
the Tracy-Widom distribution. \\
(ii) In addition, suppose that $X$ has i.i.d. standard Gaussian entries. Then $s_{\widehat{A},1}$ converges in distribution to
the Tracy-Widom distribution.
\label{theorem_TW}
\end{theorem}

\textbf{Proof:} Both claims follow from Theorem \ref{theorem_fundamental} above and Johnstone's  work on the largest eigenvalues of the Wishart and Jacobi ensembles (specifically, Theorem 1.1 in \cite{johnstone01} and Theorem 1 in \cite{johnstone08}). \hfill  $\square $

\begin{table}[htbp]
\begin{tabular}{cc}
\hline
\begin{minipage}[t]{0.45\textwidth} \centering Singular values of
$\widehat{Y}$ \\ \includegraphics[width=1\textwidth]{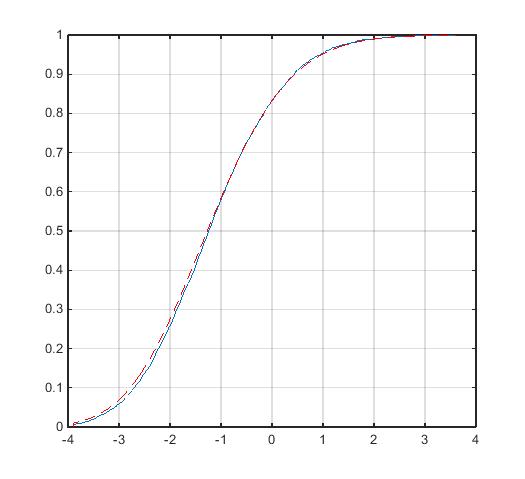}
\end{minipage} & \begin{minipage}[t]{0.46\textwidth} \centering Singular
values of $\widehat{A}$ \\
\includegraphics[width=\textwidth]{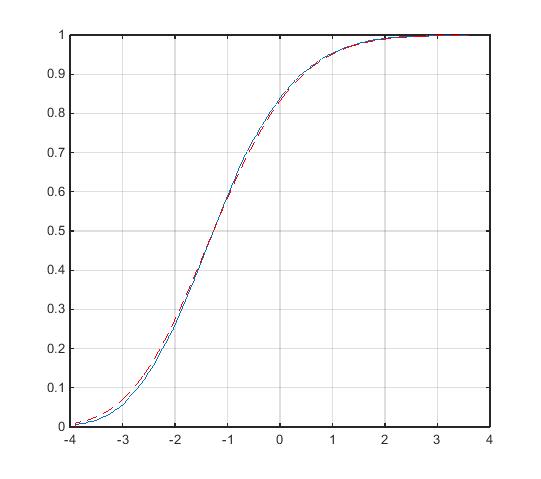} \end{minipage} \\ 
\hline
\end{tabular}%
\captionsetup{name=Figure}
\caption{ The blue solid line shows the cumulative distribution function for
the scaled largest singular value. The red dashed line is the Tracy-Widom
distribution $TW_{1}.$ The parameters are $N=100$, $p=66$, $r=133$.}
\label{fig:largest_null}
\end{table}

Corresponding to each of these three tests, we can devise a procedure for
the choice of the rank of the model (\ref{low_rank_model}):

1. Calculate the statistics $s_{Y,i}$, $s_{\widehat{Y},i}$, or $s_{\widehat{A},i}$ for the singular values of corresponding  matrices.

2. Check how many of these statistics exceed the $10\%$ quantile
of the Tracy-Widom distribution ($x_{10\%}=0.45$) and take this number as the
rank of the model (\ref{low_rank_model}).

This should be compared to the procedure suggested in \cite%
{bunea_she_wegcamp11}. The first part of their procedure is to compute $l_{\widehat{Y},i}$. Next, they prescribe to choose the
rank of the model equal to the number of statistics  $l_{\widehat{Y},i}$ that exceed a
threshold $t $. This is similar to our prescription. However, our
choice of the threshold is different from the choice in \cite%
{bunea_she_wegcamp11}. They suggest either choosing $t$ by cross-validation or using $t =2(p+r)$. We use $t =2(p+r)$
to replicate their method in numerical experiments and call this algorithm BSW (``Bunea-She-Wegkamp"). Our choice of the threshold is based on the $10$-th quantile
of the Tracy-Widom distribution.

It is not guaranteed that any of these algorithms will estimate the rank of the model successfully. First of all, as we will see later, there is a certain threshold, so that if a true singular value $\theta$ is below this threshold then it cannot be estimated consistently. Second, since Theorem \ref{theorem_TW} is about the null case, hence it is not applicable in the situation with nonzero true singular values and it does not guarantee that the largest of the remaining singular value estimates is distributed according to the Tracy-Widom law.  

The good news is that if the true singular values exceed the threshold, then Theorem 2 in \cite{onatski06} and Proposition 5.8 in \cite{BenaychGeorges_Guionnet_Maida11} suggest that our
Theorem  \ref{theorem_TW} can be extended to the non-null case by perturbative methods, and therefore the distribution of the largest of the remaining singular values of $\widehat{Y}$ or $\widehat{A}$ is indeed Tracy-Widom. (An analogous result also holds for deformations of Wigner matrices, see  Proposition 5.3 in  \cite{BenaychGeorges_Guionnet_Maida11} and Theorem 2.7 in 
\cite{knowles_yin13}.) 
 However, if one of the true singular values is at the threshold (precisely, if it is in a neigborhood of the threshold that shrinks as the matrix size grows), then the limiting distribution is a deformation of the Tracy-Widom law. 
For the real-valued Wishart matrices, the full description of this deformation is still an active research problem (\cite{forrester13a}).

In order to compare the performance of the rank-selection algorithms we run
several numerical experiments. Their results are summarized below.

In the first experiment, we assumed that the null hypothesis is satisfied
and $A=0$. The number of observations in this experiment is $N=100$, and the
number of predictors is $p=25.$

\setlength{\tabcolsep}{9pt} 
\begin{table}[htbp]
\begin{tabular}{c c |p{2cm}p{2.2cm}p{2.2cm}p{2.2cm}}
\hline
&  & BSW algorithm & Algorithm based on singular values of $Y$ & Algorithm based on singular values of $\widehat{Y}$ & Algorithm based on singular values of $\widehat{A}$ 
\\ \hline
$r=25$ &  $\overline{\widehat{s}}$ & 0.12 & 0.08 & 0.10 & 0.09 \\ 
$r=75$ &  $\overline{\widehat{s}}$ & 0.02 & 0.10 & 0.10 & 0.09 \\ \hline
\end{tabular}%
\captionsetup[table]{name=Figure}
\caption{Comparison of rank selection algorithms. \newline
Null hypothesis case: rank $s=0$. The entries in the table show the rank estimator $\widehat{s}$ averaged over $100$ repetitions of the experiment.}
\label{table_rank_selection_1}%
\end{table}

Figure \ref{table_rank_selection_1} displays the estimated rank $\widehat{s}$ averaged 
over many simulations of the model. The results show that for the null case $A=0$ the algorithms based on the Tracy-Widom
distribution falsely detect about $10\%$ of model realizations as having
rank 1 (both in the case 
$p=r=25$ and in the case $p=25<r=75$). This behavior should be expected since the threshold was set at the 10-percentile of
the Tracy-Widom distribution. In other words, the significance level of the tests is $10\%$.

 In contrast, the results of the BSW
(``Bunea-She-Wegkamp") algorithm change from about $12\%$ of false detections in
the case when $p=r=25$ to $2\%$ in the case when $p=25<r=75$. Hence, in our application of the BSW algorithm, the proportion of false positives varies depending on parameters. That is, one cannot easily pinpoint the significance level of the test based on the BSW algorithm. 

\setlength{\tabcolsep}{9pt} 
\begin{table}[htbp]
\begin{tabular}{cc|p{2cm}p{2.2cm}p{2.2cm}p{2.2cm}}
\hline
&  & BSW algorithm & Algorithm based on singular values of $Y$ & Algorithm based on singular values of $\widehat{Y}$ & Algorithm based on singular values of $\widehat{A}$ 
 \\ \hline
$r=25$ & $\overline{\widehat{s}}$ & 0.56 & 0.24 & 0.53 & 0.17 \\ 
$r=75$ & $\overline{\widehat{s}}$ & 0.78 & 0.56 & 0.92 & 0.19 \\ 
$r=200$ & $\overline{\widehat{s}}$ & 0.84 & 0.93 & 1.05 & 0.20 \\ \hline
\end{tabular}%
\captionsetup[table]{name=Figure}
\caption{Comparison of rank selection algorithms. \newline
Non-Null hypothesis case: $s=1$, $\protect\theta=0.025$}
\label{table_rank_selection_2}%
\end{table}

The results for non-null case with rank $s=1$ are shown in Figure \ref%
{table_rank_selection_2}. In this experiment, the strength of the signal $%
\theta$ was chosen to make it difficult but not impossible to detect the
signal. As before, $N=100$ and $p=25.$ 

The values in Figure \ref{table_rank_selection_2} are simulation-based estimates of the expected values of the rank estimators. For our choice of parameters, in most repetitions the rank estimator takes the values $0$ or $1$. 
Hence, the numbers in this table are approximations for the power of the test, that is, for the proportion of the 
times when the test detects the signal when it is in fact present. (Note, however, that this is only an 
approximation since the rank estimator $\widehat{s}$ can take values greater than $1$. In particular, this explain the value $1.05$ 
in the table.)

For $r=25$, the values in Figure \ref{table_rank_selection_2} show that the estimated expected values of the rank estimators based on the BSW and the fitted values algorithms are $0.56$ and $0.53$, respectively. The estimated expected values for the other two rank estimators are much smaller. Since the true rank is $1$, this finding suggests that the power of the tests  based on the BSW and the fitted values algorithms is larger than for the tests based on the other two rank estimators.  

As the number of responses, $r$, increases, it becomes easier to detect the
signal. The best performer for larger $r$ is the algorithm based on the
distribution of the largest singular value of the fitted responses. In
contrast, the algorithm based on the singular values of the coefficient matrix
estimator $\widehat{A}$ detects the signal poorly. The
performance of the BSW algorithm is also not very satisfactory as it appears
to be too conservative and biased in favor of the null hypothesis. For
example, for $r=200$, the BSW algorithm is outperformed by the simple
algorithm based on the singular values of responses, which does not use any
information about the design matrix $X$.

\setlength{\tabcolsep}{9pt} 
\begin{table}[htbp]
\begin{tabular}{cc|p{2cm}p{2.2cm}p{2.2cm}p{2.2cm}}
\hline
&  & BSW algorithm & Algorithm based on singular values of $Y$ & Algorithm based on singular values of $\widehat{Y}$ & Algorithm based on singular values of $\widehat{A}$ 
\\ \hline
$r=25$ &  $\overline{\widehat{s}}$ & 6.25 & 5.00 & 6.21 & 4.51 \\ 
$r=75$ &  $\overline{\widehat{s}}$ & 8.26 & 7.41 & 8.51 & 6.07 \\ 
$r=200$ &  $\overline{\widehat{s}}$ & 9.06 & 9.03 & 9.69 & 6.97 \\ \hline
\end{tabular}%
\captionsetup[table]{name=Figure}
\caption{Comparison of rank selection algorithms. \newline
Non-Null hypothesis case: $s=10$, $\protect\theta_1=\ldots=\protect\theta%
_{10}=0.05$;  the entries in the table show the rank estimator $\widehat{s}$ averaged over $100$ repetitions of the experiment.}
\label{table_rank_selection_3}%
\end{table}

Finally, we consider the setup, in which $A$ has $10$ nonzero singular values, each with
value $\theta = 0.05.$ The results, summarized in Figure \ref%
{table_rank_selection_3}, are similar to results for one nonzero singular value in Figure \ref{table_rank_selection_2}. The best estimates are produced by the BSW algorithm and the algorithm based on fitted
responses $\widehat{Y}$. For small $r$, the signal detection is difficult and both algorithms
perform roughly similar. For large $r$, the algorithm based on fitted
responses outperform the BSW algorithm.

\subsection{Estimation of singular values}

\subsubsection{Limit distribution of the squared singular values}

Recall that the empirical eigenvalue distribution of a square Hermitian $r$-by-$r$ matrix $M$ is the measure $\mathbb{P}=\frac{1}{r}\sum_{i=1}^{r}\delta_{  \lambda_i} $ where $\lambda _i$ are eigenvalues of $M$ and $\delta_x$ is the Dirac measure concentrated on $x$, that is, for every continuous function $f$, $\int_{\R} f \delta_x = f(x).$

\begin{table}[htbp]
\begin{tabular}{cc}
\hline
\begin{minipage}[t]{0.45\textwidth} \centering Singular values of
$\widehat{Y}$ \\
\includegraphics[width=1%
\textwidth]{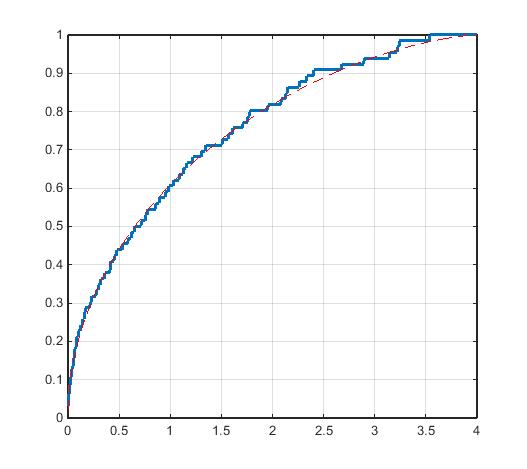} \end{minipage} & %
\begin{minipage}[t]{0.45\textwidth} \centering Singular values of
$\widehat{A}$ \\
\includegraphics[width=%
\textwidth]{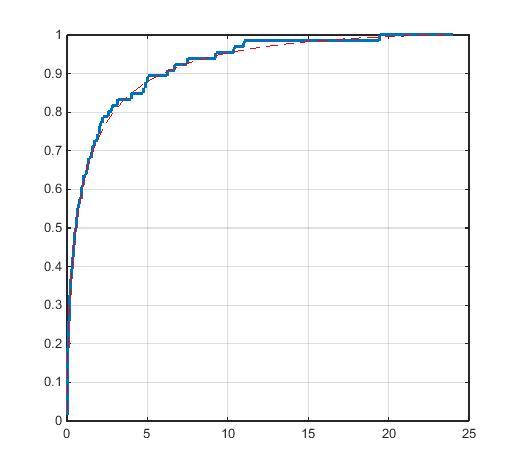} \end{minipage} \\ 
\hline
\end{tabular}%
\captionsetup{name=Figure}
\caption{ The blue line shows the cumulative distribution function of
squared singular values in a random realization of the model. The red dashed
line is the limit distribution. The parameters are $N=100$, $p=r=66$.}
\label{table_limit_distr}
\end{table}

\begin{theorem}
\label{Theorem_main}Let $\widehat{A}\df\left( X^{\ast }X\right) ^{-1}X^{\ast
}U$ and $\widehat{Y}\df X\widehat{A}.$  
\\
(i)Suppose that Assumption A1 is satisfied, $U$ is a matrix with
i.i.d. standard Gaussian entries, and $X$ is a full-rank matrix independent of $U.$ Then, as $
N\to \infty ,$ the empirical eigenvalue distribution of the $r$-by-$%
r $ matrix  $\frac{1}{r}\widehat{Y}%
^{\ast }\widehat{Y}$ converges weakly to the Marchenko-Pastur probability measure $\mathbb{P}_{MP}^{(\beta) }$ with the density defined in (\ref{density_MP}).
 \\
(ii) In addition, suppose that the entries of $X$ are i.i.d. standard Gaussian variables. Then the empirical eigenvalue distribution of $\widehat{A}^{\ast }\widehat{A}$ weakly converges to the probability measure $\mathbb{P}^{\left(\lambda ,\beta \right) }$ with the density defined in  
(\ref{density}).
\end{theorem}

The Marchenko-Pastur measure $\mathbb{P}_{MP}^{(\beta) }$ is supported on the interval
\begin{equation*}
\left[
\left( \sqrt{\beta }-1\right) ^{2},\left( \sqrt{\beta }+1\right) ^{2}\right] 
\end{equation*}
 and has the density 
\begin{equation}
p_{MP}^{\beta}(x) =\frac{1}{2\pi }\frac{1}{x}\sqrt{4\beta -\left( \beta
+1-x\right) ^{2}}. \label{density_MP}
\end{equation}%
If $0<\beta <1,$ then the measure $\mathbb{P}_{MP}^{\left( \beta \right) }$ has an
additional atom at $0$ with mass $1-\beta .$

The family of probability measures $\mathbb{P}^{\left( \lambda ,\beta \right) }$ is
parameterized by $\lambda \geq 0$ and $\beta >0.$ The continuous part of $%
\mathbb{P}^{\left( \lambda ,\beta \right) }$ is supported on the interval $\left[
x_{1},x_{2}\right] .$ If $\lambda >0,$ then 
\begin{equation}
x_{1,2}=\frac{1}{\lambda ^{2}\beta }\left[ \left( 1+\beta \right) \lambda
+2\pm 2\sqrt{\lambda ^{2}\beta +\lambda \left( \beta +1\right) +1}\right] ,
\label{support}
\end{equation}%
where $x_1$ and $x_2$ correspond to the $-$ and $+$ signs before the square root, respectively.
If $\lambda =0,$ then $x_{1}=\left( \beta -1\right) ^{2}/\left( 4\beta
\right) $ and $x_{2}=\infty .$ In both cases, the density is 
\begin{equation}
p^{(\lambda ,\beta)}(x) =\frac{1}{2\pi }\frac{1}{x\left( x+1\right) }\sqrt{%
-\left( \beta -1\right) ^{2}+2\beta \left[ \left( 1+\beta \right) \lambda +2%
\right] x-\beta ^{2}\lambda ^{2}x^{2}}.  \label{density}
\end{equation}%
If $0<\beta <1,$ then the measure $\mathbb{P}^{\left( \lambda ,\beta \right) }$ has
an additional atom at $0$ with mass $1-\beta .$

\textbf{Proof:} Both parts immediately follow from Theorem \ref{theorem_fundamental}. The first part uses a property of the Wishart random matrices discovered by Marchenko and Pastur \cite{marchenko_pastur67} 
(and independently rediscovered by Jonnson \cite{jonnson82} and Wachter \cite{wachter78}).

 The second part uses a property of the Jacobi ensemble shown by Wachter \cite{wachter80} 
and Silverstein \cite{silverstein85}. In appendix, we give another proof of this property which uses the S-transform from free probability.

\textbf{Remarks:} 1. The theorem is illustrated by a numerical example in
Figure \ref{table_limit_distr}.

2. Here is what happens in some special cases:

a) If $\lambda \rightarrow \infty $ and $\beta $ is fixed, then both $x_1$ and $x_2$ converge to $0.$ This is in agreement with the intuitive notion that
the true matrix $A_{N}=0$ can be estimated precisely if the number of
observations is large relative to the number of variables in the model.

b) If $\lambda \to \infty $ and $\beta \to 0$ so that $%
\lambda \beta \to \xi >0,$ then both $x_1$ and $x_2$ converge to $\xi^{-1} .$ Note
that $\xi =\lim_{N\to \infty }\left( N/r\right) .$

c) If $\lambda $ is fixed and $\beta \to \infty $, then both
$x_1$ and $x_2$ converge to $\lambda ^{-1}.$ 

In both b) and c),  the regression will pick up spurious dependencies. In b), it is because the number of responses is very large, and in c), it is because the number of predictors is comparable to the number of observations.

d) If $p=r$, then $\widehat{A}$ is square and one can ask about the distribution of its eigenvalues. 
By using the methods from \cite{guionnet_krishnapur_zeitouni11} and \cite{haagerup_larsen00},
 the limit distribution of eigenvalues of $\widehat{A}$ can be recovered from that of its singular values. It turns out that the limit distribution is supported
on the disc $|z|^2 \leq \frac{1}{\lambda}$ and has the density 
			\begin{equation*}
			  \frac{1+\lambda}{\pi}\frac{dm(z)}{(1+|z|^2)^2},
			\end{equation*}
where $dm(z)$ is the Lebesgue measure on the complex plane. 
After the stereographic projection this measure maps to the uniform measure on a Riemann sphere's cap.

This is a generalization of the result for the spherical ensemble of random matrices, which occurs when $N=p=r$, and therefore $\widehat{A}=X^{-1}U$. For the eigenvalues of this ensemble, it is known that the limit distribution is uniform on the Riemann sphere after the stereographic projection. (See \cite{krishnapur08}, \cite{bordenave11}, \cite{forrester_mays12} 
and \cite{tikhomirov13} for more results about this ensemble).

3. In \cite{tikhomirov13}, Tikhomirov considers the spherical ensemble $X^{-1}U$ and finds the limit
distribution of its singular values under rather weak assumptions on the
distribution of matrix entries. The basis for the results in \cite{tikhomirov13}
is the Gaussian case which is extended to non-Gaussian matrices by
consecutively changing every matrix entry to a Gaussian random variable and
then verifying that the total change in the Stieltjes transform of the
eigenvalue distribution is negligible. We conjecture that the result in
Theorem \ref{Theorem_main} can be extended to non-Gaussian matrices in a
similar way.

\subsubsection{Consistent estimation of singular values}

Now, suppose that the true matrix $A$ is from the reduced-rank model (%
\ref{low_rank_model}). We are interested to know if the singular values $%
\theta _{i}$ can be estimated consistently.

The estimator $\widehat{A}=\left( X^{\ast }X\right) ^{-1}\left( X^{\ast
}Y\right) $ is an $m$-by-$r$ matrix that can be written as follows: 
\begin{equation*}
\widehat{A}=A+\left( X^{\ast }X\right) ^{-1}\left( X^{\ast }U\right) .
\end{equation*}%
In other words, $\widehat{A}$ is the sum of the matrix $A$ and a random matrix. This random matrix has a rotationally invariant distribution by the assumption that both $X$ and $U$ are Gaussian  with i.i.d entries. In addition, matrix $A$ has a fixed set of non-zero singular values 
$\theta_1\geq \ldots \geq \theta_s > 0$ and we assume that $N\gg s$. In this situation,  one can apply the results by
Benaych-Georges and Nadakuditi from \cite{BenaychGeorges_Nadakuditi12}.

It turns out that the singular values of $\widehat{A}$ do not converge to that of $A$ as $N\to \infty$. We need to define a correction function, $D_A(x),$ that will map the singular values of $\widehat{A}$ to consistent estimates of the singular values of $A$.

 Recall that the
Stieltjes transform of a probability measure $\mathbb{P}$ is an analytic
function defined \ as 
\begin{equation*}
G\left( z\right) =\int \frac{\mathbb{P}\left( dt\right) }{z-t},
\end{equation*}%

Let $G_A(z)$ denote the  Stieltjes transform for the measure $\mathbb{P}^{(\lambda,\beta)}$
 (from Theorem \ref{Theorem_main}). It can be computed as
\begin{eqnarray}
G_A(z) &  =& \frac{1}{2z\left( 1+z\right) } \times \label{limit_Stieltjes} \\
       && \left[ 1-\beta +\left(
2+\lambda \beta \right) z-\sqrt{\left( \lambda \beta z\right) ^{2}-2\beta %
\left[ \lambda \left( 1+\beta \right) +2\right] z+\left( \beta -1\right) ^{2}%
}\right] .  \nonumber
\end{eqnarray}
(This follows from formula (\ref{Stieltjes_P_lambda_beta}) in Appendix after rescaling by $\beta^{-1}$.)

\begin{table}
\includegraphics[width=10cm]{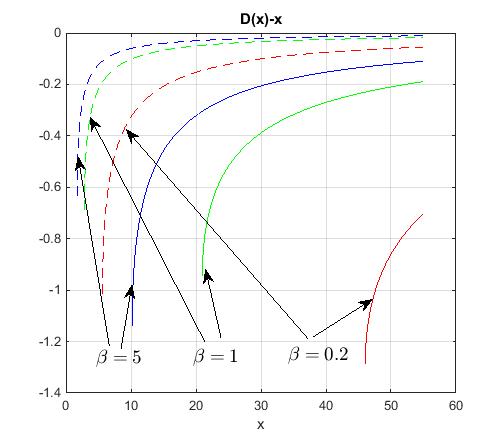} 

\captionsetup{name=Figure}
\caption{ The function $D_A(x)-x$ for various values of parameters. The dashed
lines are for $\protect\lambda=1$ and the solid lines are for $%
\protect\lambda=0.1$. The plot lines start from $x=x_2(\lambda,\beta)$ defined in (\ref{support}).
}
\label{D_function}
\end{table}

Let 
\begin{equation*}
\widetilde{G}_A(x) \df\left[ \beta ^{-1}G_A(x) +\left(
1-\beta ^{-1}\right) \frac{1}{x}\right] ,
\end{equation*}%
and 
\begin{equation*}
D_A(x) \df\frac{1}{x\sqrt{G_A\left( x^{2}\right) \widetilde{G}_A\left(
x^{2}\right) }}.
\end{equation*}%
This function has no singularities for real $x>\sqrt{x_{2}},$ (with $x_{2}$ defined
in (\ref{support})), and it is increasing on $\left(\sqrt{x_{2}},\infty \right) .$
The behavior of the function $D_A(x)-x$ for various values of parameters $%
\lambda$ and $\beta$ is illustrated in Figure \ref{D_function}.

Let 
\begin{equation*}
\overline{\theta}_A\df\lim_{x\downarrow \sqrt{x_{2}}}D_A(x).
\end{equation*}
Note that $\overline{\theta}_A\leq x_{2}$. (The threshold is below the upper edge of the support for the limit measure
 $\mathbb{P}^{(\lambda,\beta)}$.)

\begin{theorem}
\label{theorem_perturbation}
Let $\widehat{\sigma}_1 \geq \widehat{\sigma}_2 \geq \ldots \geq \widehat{\sigma}_s$  be the first $s$
largest singular values of $\widehat{A}=\left( X^{\ast }X\right) ^{-1}\left(
X^{\ast }Y\right) $, where $Y=XA+U$ and $A=\sum_{j=1}^{s}\theta _{j}u_{j}v_{j}^{\ast }$ with $\theta_1\geq \ldots \geq \theta_s > 0$. Suppose Assumptions A1 is satisfied with $\lambda
>0$ and $X$ and $U$ are independent matrices with
i.i.d. standard Gaussian entries.  For each fixed $i$, if $\theta _{i}>\overline{\theta}_A,$ then almost
surely as $N\to \infty ,$ 
\begin{equation*}
\widehat{\theta}_{i}\df D_A\left( \widehat{\sigma }_{i}\right) \to
\theta _{i}.
\end{equation*}%
If $\theta _{i}\leq \overline{\theta }_A$, then $\widehat{\sigma }%
_{i}\rightarrow \sqrt{x_{2}}.$
\end{theorem}

In other words, $\widehat{\theta }_{i}=D_A\left( \widehat{\sigma }_{i}\right) $
is a consistent estimator of $\theta _{i}$ provided that $\theta _{i}>%
\overline{\theta }_A,$ otherwise $\theta _{i}$ is hidden by spurious eigenvalues of $\widehat{A}$.

There is a similar result that uses $\widehat{Y}$ instead of $\widehat{A}$. Define 
\begin{equation*}
D_Y(x) \df\frac{1}{x\sqrt{G_{MP}^{\beta}\left( x^{2}\right) \widetilde{G}_{MP}\left(
x^{2}\right) }}, 
\end{equation*}%
where $G_{MP}^{\beta}(x)$ is the Stieltjes transform of the Marchenko-Pastur distribution with parameter $\beta$:
\begin{equation}
G_{MP}^{\beta}(z)=\frac{1}{2z}\left[ 1-\beta +z-\sqrt{(1-\beta+z)^2-4z}\right] . 
 \label{MP_Stieltjes}
\end{equation} 
and $\widetilde{G}_{MP}^{\beta}\df\left[ \beta ^{-1}G_{MP}^{\beta}(x) +\left(1-\beta ^{-1}\right) \frac{1}{x}\right].$

		   Define $x_{u,Y}$ as the upper edge of the support of the Marchenko-Pastur distribution $\mathbb{P}_{MP}^{(\beta)}$ and 
    $ \overline{\theta_Y}\df\lim_{x\downarrow \sqrt{x_{u,Y}}}D_Y(x).$ 
		
		\begin{theorem}
		\label{theorem_perturbation_Y}
		Let $\widehat{Y}=X\left( X^{\ast }X\right) ^{-1}\left(
X^{\ast }Y\right) ,$  where $Y=XA+U$ and $A=\sum_{j=1}^{s}\theta _{j}u_{j}v_{j}^{\ast }$ with $\theta_1\geq \ldots \geq \theta_s > 0$.  Suppose Assumption A1 is satisfied with $\lambda
>0$, and $X$ and $U$ are independent matrices with
i.i.d. standard Gaussian entries. Let $\widehat{\sigma }_{i,\widehat{Y}}$ be the first $s$ largest singular values of $\widehat{Y}/\sqrt{r}.$
			If $\theta _{i}>\overline{\theta_Y},$ then almost surely as $N\rightarrow \infty ,$ 
      \begin{equation*}
        \widehat{\theta }_{i,\widehat{Y}}\df \mu^{-1/2}D_Y\left( \widehat{\sigma }_{i,\widehat{Y}}\right) \to
        \theta _{i}.
      \end{equation*}%
			
      If $\theta _{i}\leq \overline{\theta_Y }$, then $\widehat{\sigma}_{i,\widehat{Y}}\rightarrow \sqrt{x_{u,Y}}.$
		\end{theorem}

\begin{table}[tbph]
\includegraphics[width=9cm]{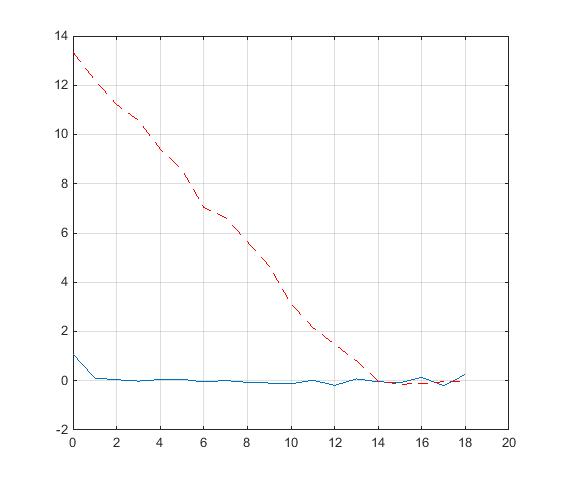}
\captionsetup{name=Figure}
\caption{ Quality of the estimators for the rank-one perturbation. The
horizontal axis shows the parameter $\protect\theta$. The red dashed line is
the difference $\widehat{\protect\theta}_A - \protect\theta$ and the blue
solid line is the difference $\widehat{\protect\theta}_Y-\protect\theta$.
The thresholds are $\overline{\protect\theta}_A=14.2$ and $\overline{\protect\theta}_Y=0.8$. (Parameters are $%
N=500$, $\protect\lambda=0.2$, $\protect\beta=1/2$.)}
 \label{fig_estimate}
\end{table}		

\begin{table}[htbp]
  \begin{tabular}{cc}
    \hline
    \begin{minipage}[t]{0.5\textwidth} 
		  \includegraphics[width=\textwidth]{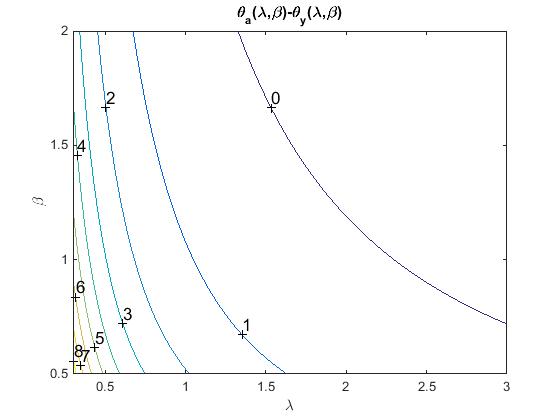} 
		\end{minipage} & %
    \begin{minipage}[t]{0.5\textwidth} 
      \includegraphics[width=\textwidth]{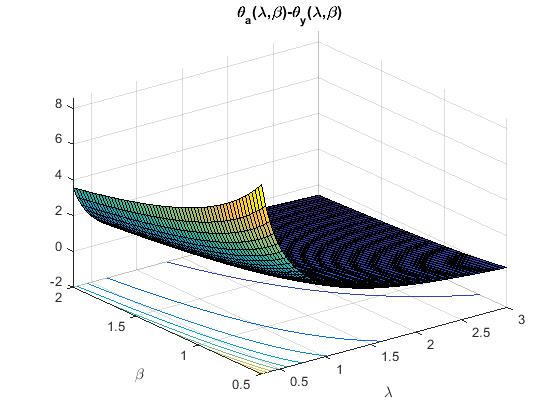} 
	  \end{minipage} \\ 
    \hline
  \end{tabular}%
	\captionsetup{name=Figure} 
  \caption{ Contour and surface plots of the difference between the estimator thresholds 
	$\overline{\theta}_A(\lambda,\beta)$ and $\overline{\theta}_Y(\lambda,\beta)$ .}
  \label{table_thresholds}
\end{table}

\textbf{Remarks}: 1. These results are similar in spirit to the results in \cite%
{baik_benarous_peche05}, \cite{baik_silverstein06}, and \cite{paul07}, which
are concerned with the singular values of sample covariance matrices. In
these papers, it was found that if the true covariance matrix has a 'spike'
that exceeds a certain threshold, then it will be observed in the spectrum
of the sample covariance matrix. Otherwise it will be hidden among the
spurious eigenvalues.

2. The quality of the estimators is illustrated in Figures \ref{fig_estimate} and \ref{table_thresholds}. They show that for relatively small values of parameters $\lambda$ and $\beta$, the threshold $\overline{\theta}_Y$ is smaller than $\overline{\theta}_A$ and the estimator  $\widehat{\theta }_{i,\widehat{Y}}$ is preferable to $\widehat{\theta }_{i,\widehat{A}}$. For large $\lambda$ and $\beta$ the threshold $\overline{\theta}_A$ can be smaller than $\overline{\theta}_Y$. However, the difference is small and in this region the correction term $D(x)-x$ is also small.  

3. The proofs of Theorem \ref{theorem_perturbation} and Theorem \ref{theorem_perturbation_Y} are essentially by combining Theorem %
\ref{Theorem_main} in this paper and Theorem 2.8 in \cite%
{BenaychGeorges_Nadakuditi12}. We provide a detailed proof in Section \ref%
{subsection_largest_eigenvalue_nonnull} below for the convenience of the
reader. (The proof of Theorem 2.8 in \cite{BenaychGeorges_Nadakuditi12}
lacks some details and has annoying typos.)
However, we defer to \cite{BenaychGeorges_Nadakuditi12} for details about convergence and continuity issues.

\subsubsection{CLT for the estimator of the singular values}

 Let us define the empirical version of the function $G_A(x)$: 
\begin{equation}
G_{A,N}\left( x\right) \df\frac{1}{p}\sum\limits_{i=1}^{p}\frac{1}{x-s_{i}^{2}},
\label{def_G_N}
\end{equation}%
where $s_{i}$ are singular values of the matrix $X\backslash U.$ Let 
\begin{equation}
\widetilde{G}_{A,N}(x) \df\frac{r}{p}G_{A,N}(x) +\left( 1-\frac{r}{p}\right) \frac{1}{x}  \label{def_G_N_tilda}
\end{equation}%
and
\begin{equation}
D_{A,N}\left( x\right) \df\frac{1}{x\sqrt{G_{A,N}\left( x^{2}\right) \widetilde{G}%
_{A,N}\left( x^{2}\right) }}.  \label{def_D_N}
\end{equation}
Then Theorem \ref{theorem_perturbation} holds with $D_{A,N}\left( x\right) $
instead of $D_A\left( x\right) .$

Next, let $\widehat{\theta }^{\left( N\right) }\df D_{A,N}\left( \widehat{\sigma }%
^{\left( N\right) }\right) ,$ where $\widehat{\sigma }^{\left( N\right) }$
is the largest singular value of the matrix $\widehat{A}_{N}=\left(
X_{N}^{\ast }X_{N}\right) ^{-1}\left( X_{N}^{\ast }Y_{N}\right) $ and $D_{A,N}$
is as defined in (\ref{def_D_N}).

\begin{theorem}
\label{theorem_clt}
Assume  that $Y=XA+U$ with $A=\theta u v^{\ast }$. Suppose Assumption A1 is satisfied with $\lambda
>0$, and $X$ and $U$ are independent matrices with
i.i.d. standard Gaussian entries.

Then 
\begin{equation*}
\frac{\sqrt{r}}{\omega}\left( \widehat{\theta }^{\left( N\right) }-\theta
\right),
\end{equation*}%
with $\omega$ as defined below in (\ref{omega}), converges in distribution to a standard zero-mean Gaussian random variable.
\end{theorem}

\begin{table}[tbph]
\includegraphics[width=10cm]{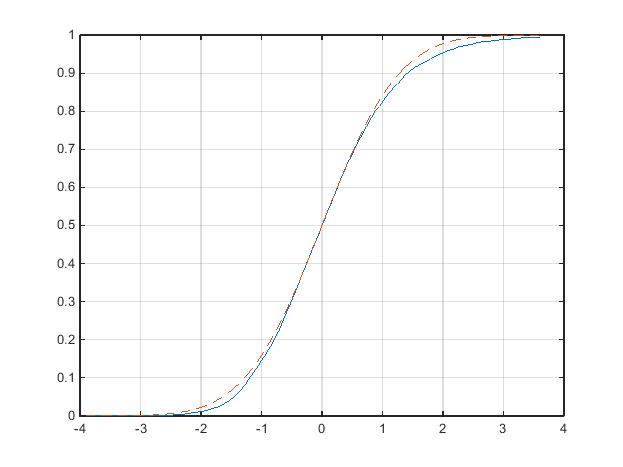} 
\captionsetup{name=Figure}
\caption{The blue solid line shows the cumulative distribution function for $%
\frac{\protect\sqrt{r}}{\protect\omega}(D(\widehat{\protect\sigma}^{(N)})-%
\protect\theta).$ The red dashed line is the standard Gaussian cumulative
distribution function. The parameters are $N=300$, $\protect\lambda=0.2$, $%
\protect\beta=2$, and $\protect\theta=9.$}
\label{fig:largest_non-null}
\end{table}

Let $\sigma $ be the largest solution of the
equation 
\begin{equation*}
D_A(x) =\theta .
\end{equation*}%
Define 
\begin{eqnarray*}
\kappa _{1}^{2} &=&-2\sigma ^{2}\left( G_A^{\prime }\left( \sigma ^{2}\right)
+ \left[ G_A\left( \sigma ^{2}\right) \right] ^{2}\right) , \\
\kappa _{2}^{2} &=&-\frac{2}{\beta }\sigma ^{2}\left( \widetilde{G}_A^{\prime
}\left( \sigma ^{2}\right) +\left[ \widetilde{G}_A\left( \sigma ^{2}\right) %
\right] ^{2}\right) , \\
\tau ^{2} &=&-\left[ \sigma ^{2}\widetilde{G}_A^{\prime }\left( \sigma
^{2}\right) +\widetilde{G}_A\left( \sigma ^{2}\right) \right].
\end{eqnarray*}%

Let 
\begin{eqnarray}
\kappa ^{2} &=&\sigma ^{2}\left[ \widetilde{G}_A \left( \sigma ^{2}\right) \right]
^{2}\kappa _{1}^{2}+\sigma ^{2}\left[ G_A\left( \sigma ^{2}\right) %
\right] ^{2}\kappa _{2}^{2}+4\theta ^{-2}\tau ^{2},\text{ and} \notag \\
\omega &=&\frac{1}{2}\theta ^{3}\kappa . \label{omega}
\end{eqnarray}

Remark: We conjecture that the statement of the theorem holds true with the estimator $%
D\left( \widehat{\sigma }^{\left( N\right) }\right) $ instead of $%
D_{N}\left( \widehat{\sigma }^{\left( N\right) }\right) .$ However the proof
of this statement runs into some technical difficulties and might require
some additional assumption on the convergence $N/m\rightarrow 1+\lambda $
and $m/r\rightarrow \beta .$

\section{Proofs of the main results}

\label{subsection_largest_eigenvalue_nonnull}

\subsection{Proof of Theorem \ref{theorem_perturbation}}
 
The basic tool is the following determinantal identity. Let $|X|$ denote the determinant of a matrix $X$.

\begin{lemma}
\label{lemma_perturbation1}Let $\Psi$ and $D$ be an $N$-by-$N$ and $s$-by-$s$
matrices, respectively, and let $W_{1}$ and $W_{2}$ be an $N$-by-$s$ and $s$-by-$N$ matrices, respectively. Assume that $D$ is invertible. Then,

\begin{equation*}
 \left| tI_N-\Psi-W_{1}DW_{2}\right| =|D|\left| tI_N-\Psi\right|
  \left| D^{-1}-W_{2}\left( tI_N-\Psi\right) ^{-1}W_{1}\right| .
\end{equation*}
\end{lemma}

\textbf{Proof: }The proof is through a sequence of elementary determinantal
identities: 
\begin{eqnarray*}
 \left| tI_{N}-\Psi-W_{1}DW_{2}\right| &=&\left| tI_{N}-\Psi\right|
\left| I_{N}-\left( tI_{N}-\Psi\right) ^{-1}W_{1}DW_{2}\right| \\
&=& \left| tI_{N}-\Psi\right| \left| I_{r}-DW_{2}\left(
tI_{N}-\Psi\right) ^{-1}W_{1}\right| \\
&=& \left| tI_{N}-\Psi\right| |D| \left| D^{-1}-W_{2}\left(
tI_{N}-\Psi\right) ^{-1}W_{1}\right| .
\end{eqnarray*}%
\hfill $\square $

We apply this lemma to understand how the singular values of a matrix $Z$ are affected if $Z$ is perturbed by a low-rank matrix $A$.

\begin{lemma}
\label{lemma_reduction_to_M}Let $\widehat{A}=A+Z$ with $A=\sum_{j=1}^{s}\theta _{j}u_{j}v_{j}^{\ast }.$ Let $\Theta =\mathrm{diag}\left( \theta _{1},\ldots ,\theta _{s}\right) ,$
and let $U$ and $V$ be an $p$-by-$s$ and an $r$-by-$s$ matrices whose
columns are vectors $u_{i}$ and $v_{i}$,
respectively. Let 
\begin{eqnarray}
M\left( t\right) &\df &\left( 
\begin{array}{cc}
tV^{\ast }\left( t^{2}I_{r}-Z^{\ast }Z\right) ^{-1}V, & V^{\ast }Z^{\ast
}\left( t^{2}I_{p}-ZZ^{\ast }\right) ^{-1}U \notag \\ 
U^{\ast }\left( t^{2}I_{p}-ZZ^{\ast }\right) ^{-1}ZV, & tU^{\ast }\left(
t^{2}I_{p}-ZZ^{\ast }\right) ^{-1}U%
\end{array}%
\right) \\
&-&\left( 
\begin{array}{cc}
0 & \Theta ^{-1} \\ 
\Theta ^{-1} & 0%
\end{array}%
\right) .  \label{definition_Mt}
\end{eqnarray}
 If $t>0$ and $\left|
M\left( t\right) \right|=0,$ then $t$ is a singular value of $\widehat{A}.$
Conversely, if a singular value of $\widehat{A}$ is different from the
singular values of $Z$, then it is a zero of the function $\left| M\left( t\right)\right|$ defined in  (\ref{definition_Mt}) below.
\end{lemma}

\textbf{Proof:} We apply Lemma (\ref{lemma_perturbation1}) to $\Psi=\left( 
\begin{array}{cc}
0 & Z \\ 
Z^{\ast } & 0%
\end{array}%
\right) ,$ $D=\left( 
\begin{array}{cc}
0 & \Theta \\ 
\Theta & 0%
\end{array}%
\right) ,$ $W_{1}=\left( 
\begin{array}{cc}
0 & U \\ 
V & 0%
\end{array}%
\right) ,$ and $W_{2}=W_{1}^{\ast },$ and we use the identity: 
\begin{equation*}
\left( 
\begin{array}{cc}
tI_{p} & -Z \\ 
-Z^{\ast } & tI_{r}%
\end{array}%
\right) ^{-1}=\left( 
\begin{array}{cc}
t\left( t^{2}I_{p}-ZZ^{\ast }\right) ^{-1} & \left( t^{2}I_{p}-ZZ^{\ast
}\right) ^{-1}Z \\ 
Z^{\ast }\left( t^{2}I_{p}-ZZ^{\ast }\right) ^{-1} & t\left(
t^{2}I_{r}-Z^{\ast }Z\right) ^{-1}%
\end{array}%
\right) .
\end{equation*}
Then we get 
\begin{eqnarray*}
 \left| tI_{p+r}-\left( 
\begin{array}{cc}
0 & A+Z \\ 
A^{\ast }+Z^{\ast } & 0%
\end{array}%
\right) \right| &=&\left( -1\right) ^{s}\prod\limits_{i=1}^{s}\theta
_{i}^{2} \left| tI_{p+r}-\left( 
\begin{array}{cc}
0 & Z \\ 
Z^{\ast } & 0%
\end{array}%
\right) \right| 
\\
&&\times  \left|M\left( t\right)\right| .
\end{eqnarray*}%
Hence, the positive eigenvalues of the matrix $\left( 
\begin{array}{cc}
0 & A+Z \\ 
A^{\ast }+Z^{\ast } & 0%
\end{array}%
\right)$ are either positive eigenvalues of the matrix $\left( 
\begin{array}{cc}
0 & Z \\ 
Z^{\ast } & 0%
\end{array}%
\right) $ or positive zeros of $\left|M\left( t\right)\right| .$ Recall that the
singular values of matrix $A+Z$ coincide with the positive eigenvalues of
matrix 
\begin{equation*}
\left( 
\begin{array}{cc}
0 & A+Z \\ 
(A+Z)^{\ast } & 0%
\end{array}%
\right) ,
\end{equation*}
and similarly for the singular values of matrix $Z.$ This proves both claims
of the Lemma. \hfill $\square $

We will apply this theorem to $Z=X\backslash \mathcal{U}$, where $\mathcal{U}$ denotes the matrix of noise in the regression model \ref{regression_model} and  should not be confused with matrix $U$ in the definition of $M(t)$. The crucial observation is that the rotational invariance of $Z$ and Theorem \ref{Theorem_main} imply that the $s$-by-$s$ matrix $%
tV^{\ast }\left( t^{2}I_{r}-Z^{\ast }Z\right) ^{-1}V\ $ converges to the
scalar matrix $tG_A\left( t^{2}\right) I_{s}$ provided that $t>\sqrt{x_2}$. 

By using Lemma \ref{lemma_shift}, one can also check that the matrix \\
 $tU^{\ast }\left( t^{2}I_{p}-ZZ^{\ast }\right) ^{-1}U$ converges to $\left[
\beta ^{-1}tG_A\left( t^{2}\right) +\left( 1-\beta ^{-1}\right) t^{-1}\right]
I_{s}.$

The off-diagonal blocks of $M\left( t\right) $ converge to zero.

Finally, we use the fact that 
\begin{equation*}
\left| 
\begin{array}{cc}
xI_{s} & -\Theta ^{-1} \\ 
-\Theta ^{-1} & yI_{s}%
\end{array}%
\right| =\prod\limits_{i=1}^{s}\left( xy-\theta _{i}^{-2}\right) .
\end{equation*}%
and conclude (by continuity properties of the determinantal equation
solutions proved in \cite{BenaychGeorges_Nadakuditi12}) that if $\theta_i>\overline{\theta}_A$ then there is a singular value of $\widehat{A}$ that converge to the positive solution of the
following equation:%
\begin{equation*}
tG_A\left( t^{2}\right) \left[ \beta ^{-1}tG_A\left( t^{2}\right) +\left(
1-\beta ^{-1}\right) t^{-1}\right] =\theta _{i}^{-2}.
\end{equation*}%
\hfill Since the left hand side equals $\left[ D_A\left( t\right) \right]
^{-2}$ by definition, this convergence establishes the first part of Theorem \ref%
{theorem_perturbation}. 

In order to establish the second part, let 
 $\theta_1  \geq \ldots \geq \theta_k > \overline{\theta_A} \geq \theta_{k+1} \geq \ldots \geq \theta_s$. Then, by monotonicity of $D_A(x)$, if $i>k$ then for every $\varepsilon>0$ and sufficiently large $N$ the equation $D_A(x)=\theta_i$ does not have roots in the interval $(\sqrt{x_2}+\varepsilon,\infty)$. Hence, by Lemma \ref{lemma_reduction_to_M}), the eigenvalues of matrix $\widehat{A}$ can have only $k$ limits greater than $\sqrt{x_2}+\varepsilon$. It can also be checked that the singular values of the matrix $\widehat{A}$ has the same limiting distribution, $\bP^{(\lambda,\beta)}$, as in the null case. (Indeed, a deformation by a fixed-rank matrix does not affect the empirical distribution of eigenvalues). It follows that the $(k+1)$-st, $(k+2)$-nd, \ldots, and $s$-th largest singular values of $\widehat{A}$ must converge to  $\sqrt{x_2}$.
         \hfill $\square $

\subsection{Proof of Theorem \ref{theorem_perturbation_Y}} 

An analogue of Lemma \ref{lemma_reduction_to_M} holds and establishes the fact that the singular values of $\widehat{Y}/\sqrt{r}$ are zeros of $\left|M_Y(t)\right|$, where
\begin{eqnarray*}
M_Y(t) &=& -\left( 
\begin{array}{cc}
0 & \Theta ^{-1} \\ 
\Theta ^{-1} & 0%
\end{array}%
\right)\\
&+&\left( 
\begin{array}{cc}
  tV^{\ast }\left( t^{2}I_{r}-Z^{\ast }Z\right) ^{-1}V, &
  r^{-\frac{1}{2}}V^{\ast }Z^{\ast}\left( t^{2}I_{N}-ZZ^{\ast }\right) ^{-1}XU \\ 
  r^{-\frac{1}{2}}U^{\ast }X^{\ast}\left( t^{2}I_{N}-ZZ^{\ast }\right) ^{-1}ZV, & 
	r^{-1}tU^{\ast }X^{\ast}\left(t^{2}I_{N}-ZZ^{\ast }\right) ^{-1}XU 
\end{array}%
\right),  \label{definition_M_Y_t}
\end{eqnarray*}
where $Z= r^{-\frac{1}{2}}X(X^{\ast}X)^{-1}X \mathcal{U}$, and $\mathcal{U}$ is the matrix of noise in the regression model \ref{regression_model}.
As before, the off-diagonal blocks of $M_Y(t)$ converge to $-\Theta ^{-1}$. Besides, if $t$ is sufficiently large, then the first diagonal block converges to $tG_{MP}^{\beta}(t^2)I_s.$ 

In addition, for sufficiently large $t$, the limit of the second diagonal block is the scalar matrix $a t I_s$, where $a$ is the limit of 
\begin{equation}
\frac{N}{r}\frac{1}{p} \mathrm{Tr} \left( N^{-1}X^{\ast}\left(t^{2}I_{N}-ZZ^{\ast }\right) ^{-1}X \right) \label{matrix_trace}
\end{equation}

In order to calculate this limit, we note that the matrices
$\left(\mathcal{U}/\sqrt{r}, \mathcal{U}^{\ast}/\sqrt{r}\right)$ and 
$\left(X/\sqrt{N}, X^{\ast}/\sqrt{N}\right)$ converge in distribution (in the sense of free probability theory)
 to pairs of non-commutative random variables
 $\left(u,u^{\ast}\right)$ and $\left(x,x^{\ast}\right)$, which are free from each other. Let $P_x=x (x^{\ast}x)^{-1} x^{\ast}$. This is the orthogonal projection corresponding to $xx^{\ast}$, since $P_x xx^{\ast}=xx^{\ast}P_x=xx^{\ast}$, $P_x^2=P_x$ and $P_x^{\ast}=P_x$. Then, the limit of the expectation of (\ref{matrix_trace}) is equal to the limit of 
\begin{equation}
\frac{N}{r}(1+\lambda) \tau \left( \left(t^{2}I-P_x u u^{\ast} P_x\right) ^{-1}x x^{\ast}\right),
\label{R_XX}
\end{equation}
where $\tau$ is the trace in the corresponding free probability space.

In order to handle this expression, we use the following lemma.
\begin{lemma}
Suppose that the pair of variables $\left( a, P_a \right)$ and the variable $b$ are free from each other, and suppose that $P_a$ has the properties $P_a a= a P_a =a$ and $P_a^2=P_a$. Then 
\begin{equation*}
\tau \left((tI-P_a b P_a)^{-1} a\right)=\frac{\tau(a)}{\tau(P_a)}
\left[-\frac{1-\tau(P_a)}{t}+\tau \left((tI-P_a b P_a)^{-1}\right)\right]
\end{equation*}
\end{lemma}

\textbf{Proof:} On both sides we have complex-analytic functions in $t$ (aside of singularities), and we can expand these functions in powers of $t^{-1}$. It is enough to  show the equality in the region of the complex plane where the resulting series converge, since for other values of $t$ the equality will hold by analytic continuation. Consequently, it is enough to show that the equality holds term by term in the series. In particular, we need to check that 
\begin{equation*}
\tau \left((P_a b P_a)^n a\right)=\tau \left((P_a b P_a)^n\right) \frac{\tau(a)}{\tau(P_a)}
\end{equation*} 
for each integer $n \geq 1$.
By using the properties of $P_a$ and $a$, the left hand-side can be written as 
\begin{equation*}
\tau \left(\underbrace{P_a b \ldots P_a b}_{n\ \mathrm{times}}a\right)=\tau \left( b P_a b \ldots P_a ba\right).
\end{equation*}
Next we use the property that $\left( a, P_a \right)$ and $b$ are free, and write:
\begin{equation}
\tau \left( b P_a b \ldots P_a ba\right)=\sum_{\pi \in \mathcal{NC}(n)} \kappa_{\pi}(b,b,\ldots,b) \tau_{K(\pi)}(P_a,\ldots,P_a,a),
\label{Nica_Speicher}
\end{equation}
where the sum is extended over all non-crossing partitions of the sequence $\left\{1,2,\ldots,n\right\}$, $K(\pi)$ is the Kreweras complement of the partition $\pi$, $\tau_{K(\pi)}$ is the multiplicative extension of the trace $\tau$, associated with partition
 $K(\pi)$, and $\kappa_{\pi}$ denotes the free cumulant functional associated with partition $\pi$. (For additional details and a proof of formula (\ref{Nica_Speicher}), see Theorem 14.4 in \cite{nica_speicher06}.)

Let $|\pi|$ denote the number of blocks in the partition $\pi$. It is a fact that $|K(\pi)|=n+1-|\pi|$.
 Then, by multiplicativity of $\tau_{K(\pi)}$, the expression in (\ref{Nica_Speicher}) can be written as 
\begin{equation*}
\sum_{\pi \in \mathcal{NC}(n)} \kappa_{\pi}(b,b,\ldots,b) \tau(P_a)^{n-|\pi|} \tau(a).
\end{equation*}
We can take $\tau(a)$ outside of the sum sign, and a similar calculation shows that 
\begin{equation*}
\tau \left((P_a b P_a)^n\right)=\tau \left(P_a\right)\sum_{\pi \in \mathcal{NC}(n)} \kappa_{\pi}(b,b,\ldots,b) \tau(P_a)^{n-|\pi|}.
\end{equation*}
Hence, 
\begin{equation*}
\tau \left((P_a b P_a)^n a\right)=\tau \left((P_a b P_a)^n\right)\tau \left(a\right) /\tau \left(P_a\right).
\end{equation*}
\hfill $\square$

By applying this lemma to the expression in (\ref{R_XX}), we calculate that its limit equals  
\begin{equation*}
\mu\widetilde{G}_{MP}^{\beta}\left( t^2 \right), 
\end{equation*}
where $ \mu=\lim_{N\to\infty}\frac{N}{r}$.
 
Hence, for sufficiently large $t$, $\left|M_Y(t)\right|$ converges to 
\begin{equation*}
\prod_{i=1}^{s}\left(\mu t {G}_{MP}^{\beta}\left( t^2 \right) \widetilde{G}_{MP}^{\beta}\left( t^2 \right) - \theta_i^{-2}\right).
\end{equation*}
From this, we conclude that for $\theta_i>\overline{\theta}_Y$, there exists a singular value of $\widehat{Y}/\sqrt{r}$,
 denoted by $\widehat{\sigma}_i$, such that 
\begin{equation*}
\mu^{-1/2} D_Y(\widehat{\sigma}_i) \to \theta_i
\end{equation*}
in probability. This is equivalent to the claim in the first part of the Theorem. The second part is proved as in Theorem \ref{theorem_perturbation}. \hfill $\square$ 

\subsection{Proof of Theorem \ref{theorem_clt}} 

Now we are going to prove the
central limit theorem for the estimator $\widehat{\theta }=D\left( \widehat{%
\sigma }\right) .$ Recall that we are dealing with the rank-one perturbation
model, in which $\theta $ is the singular value of the perturbation and $%
\widehat{\sigma }$ is the largest singular value of the estimator
 $\widehat{A} =X\backslash Y.$ First, we note that the function $\frac{1}{r}\sum_{i=1}^{r}%
\frac{t}{t^{2}-s_{i}^{2}}$ concentrates better than what could be expected
if $s_{i}$ were independent.

\begin{lemma}
\label{speed_convergence}Let $s_{i},$ $i=1,\ldots ,r$ \ be the singular
values of the $m$-by-$r$ matrix $Z=X\backslash U $ and $x_{2}$ be as defined in (\ref{support}). Then for every $\varepsilon >0$ and
every $t>x_{2},$ 
\begin{equation*}
r^{1-\varepsilon }\left[ \frac{1}{r}\sum_{i=1}^{r}\frac{t}{t^{2}-s_{i}^{2}}-%
\mathbb{E}\frac{1}{r}\sum_{i=1}^{r}\frac{t}{t^{2}-s_{i}^{2}}\right]
\rightarrow 0
\end{equation*}%
in $L^{2}$.
\end{lemma}

\textbf{Proof:} The claim of this lemma is a consequence of a CLT for linear
eigenvalue statistics. Indeed, the eigenvalues of $Z^{\ast }Z$ are the
transformed eigenvalues from the Jacobi ensemble of random matrices (with a
smooth transformation (\ref{relation_of_eigenvalues})). Hence we can
use the results about the linear eigenvalue statistics of the Jacobi
ensemble. For example, the results of Dumitriu and
Paquette \cite{Dumitriu_Paquette12} (specifically their Theorem 3.1) imply
that $\mathrm{Tr}\left( t^{2}I-Z^{\ast }Z\right) ^{-1}$ converges to a
Gaussian random variable with a finite variance. Hence, $r^{-\varepsilon
}\left( \mathrm{Tr}\left( t^{2}I-Z^{\ast }Z\right) ^{-1}-\mathbb{E}\mathrm{Tr%
}\left( t^{2}I-Z^{\ast }Z\right) ^{-1}\right) \rightarrow 0$ in $L^{2}.$
\hfill $\square $

\textbf{Remark:} The results of Johansson in \cite{Johansson98} (e.g., his
Theorem 2.4) suggest that one can use $tG_A\left( t^{2}\right) $ instead of $%
\mathbb{E}\frac{1}{r}\sum_{i=1}^{r}\frac{t}{t^{2}-s_{i}^{2}}$ in the above
lemma. However, Johansson's results are proved under assumption of a fixed
potential $V,$ while the potential in the pdf of the Jacobi ensemble is changing
with $N.$ Hence, some additional assumptions on the convergence $%
m/N\rightarrow \lambda $ and $m/r\rightarrow \beta $ might be needed.

Next, we prove a CLT for the matrix $M\left( t\right) .$

\begin{lemma}
\label{lemma_M_fluctuations}Assume  that $Y=XA+U$ with $A=\theta u v^{\ast }$. Let $M_{N}\left( t\right) $
be the $2$-by-$2$ matrix defined in (\ref{definition_Mt}), and let $%
G_{N}\left( t\right) $ and $\widetilde{G}_{N}\left( t\right) $ be as defined
in (\ref{def_G_N}) and ($t$\ref{def_G_N_tilda}), respectively. If $t>x_{2}$
(with $x_{2}$ defined in (\ref{support})), then the random matrix 
\begin{equation*}
\widehat{M}_{N}\left( t\right) =\sqrt{r}\left[ M_{N}\left( t\right) -\left( 
\begin{array}{cc}
tG_{N}\left( t^{2}\right) & -\theta ^{-1} \\ 
-\theta ^{-1} & t\widetilde{G}_{N}\left( t^{2}\right)%
\end{array}%
\right) \right]
\end{equation*}%
converges in distribution to the matrix 
\begin{equation*}
\left( 
\begin{array}{cc}
\kappa _{1}X & \tau Z \\ 
\tau Z & \kappa _{2} Y%
\end{array}%
\right) ,
\end{equation*}%
where $X,Y,$ and $Z$ are independent standard Gaussian random variables, and 
\begin{eqnarray*}
\kappa _{1}^{2} &=&-2t^{2}\left( G^{\prime }\left( t^{2}\right) +\left[
G\left( t^{2}\right) \right] ^{2}\right) , \\
\kappa _{2}^{2} &=&-2\beta ^{-1}t^{2}\left( \widetilde{G}^{\prime }\left(
t^{2}\right) +\left[ \widetilde{G}\left( t^{2}\right) \right] ^{2}\right) ,
\\
\tau ^{2} &=&-\left[ t^{2}\widetilde{G}^{\prime }\left( t^{2}\right) +%
\widetilde{G}\left( t^{2}\right) \right] .
\end{eqnarray*}
\end{lemma}

For the proof we rely on Theorem 7.1 in \cite{bai_yao08} that allows one to
compute the distributional limit for the forms $\sum_{i,j=1}^{n}u_{i}\left(
l\right) A_{ij}u_{j}\left( l^{\prime }\right) ,$ as $n\rightarrow \infty ,$
for independent, identically distributed $K$-tuples of real or complex
valued random variables $\left\{ u_{i}\left( 1\right) ,\ldots ,u_{i}\left(
K\right) \right\} ,$ under some additional assumptions on matrices $\left(
A_{ij}\right) $. The theorem shows that the limit
is Gaussian and provides a formula for the covariance matrix of the limit.
This theorem is not directly applicable in our case since we will have
variables $u_{i}$ which are the coordinates of the vector uniformly
distributed on a sphere $S^{n}$ and, therefore, not independent. This can be
overcome either by a suitable modification of Theorem 7.1 in \cite{bai_yao08}%
, or by a trick that represents the uniformly distributed vector $u$ as a
normalization of a Gaussian vector (which is similar to what is done in the
proof of Theorem 6.4 in \cite{BenaychGeorges_Guionnet_Maida11}). In the
following proof we concentrate on explaining how the variance coefficients
are calculated.

\textbf{Proof of Lemma \ref{lemma_M_fluctuations}:} By a suitable
modification of Theorem 7.1 in \cite{bai_yao08}, we find that the matrix $%
\widehat{M}_{N}\left( t\right) $ converges in distribution to a matrix that
consists of independent zero-mean Gaussian entries. It remains to compute
the variance of the entries. Let us start with
the upper-left diagonal entry of $M_{N}\left( t\right) ,$ which is $tv^{\ast
}\left( t^{2}I_{r}-Z^{\ast }Z\right) ^{-1}v.$ Because of the rotational invariance of $Z,$ one can
take $v$ uniformly distributed on the sphere $S^{r}$. In the basis that
diagonalizes $Z^{\ast }Z,$ we get $\,$%
\begin{equation*}
\left[ M_{N}\left( t\right) \right] _{11}=\sum_{i=1}^{r}\frac{t}{%
t^{2}-s_{i}^{2}}v_{i}^{2},
\end{equation*}%
where $s_{i}\,$are singular values of $Z$ and $v_{i}$ are coordinates of
vector $v.$ Hence,
\begin{eqnarray*}
\mathbb{E}\left( \left[ M_{N}\left( t\right) \right] _{11}|Z\right) &=&\frac{%
1}{r}\sum_{i=1}^{r}\frac{t}{t^{2}-s_{i}^{2}},\text{ and } \\
\mathbb{E}\left( \left[ M_{N}\left( t\right) \right] _{11}\right) &=&\mathbb{%
E}\frac{1}{r}\sum_{i=1}^{r}\frac{t}{t^{2}-s_{i}^{2}}.
\end{eqnarray*}%
Next we use the total variance formula: 
\begin{equation*}
\mathbb{V}ar\left( \sqrt{r}\left[ M_{N}\left( t\right) \right] _{11}\right) =%
\mathbb{E}\left( \mathbb{V}ar\left( \left[ \sqrt{r}M_{N}\left( t\right) %
\right] _{11}|Z\right) \right) +\mathbb{V}ar\left( \mathbb{E}\left( \sqrt{r}%
\left[ M_{N}\left( t\right) \right] _{11}|Z\right) \right) .
\end{equation*}%
The second term converges to zero by Lemma \ref{speed_convergence}, hence we
need only to compute the limit of the first term. We calculate: 
\begin{equation*}
\mathbb{E}\left( r\left[ M_{N}\left( t\right) \right] _{11}^{2}|Z\right) =%
\frac{3}{(r+2)}\sum_{i=1}^{r}\left( \frac{t}{t^{2}-s_{i}^{2}}\right) ^{2}+%
\frac{1}{(r+2)}\sum_{i\neq j}^{r}\left( \frac{t}{t^{2}-s_{i}^{2}}\right)
\left( \frac{t}{t^{2}-s_{j}^{2}}\right) ,
\end{equation*}%
where we used $E\left( v_{i}^{4}\right) =\frac{3}{r(r+2)}$ and $E\left(
u_{i}^{2}u_{j}^{2}\right) =\frac{1}{r(r+2)}.$ Hence, 
\begin{eqnarray*}
\mathbb{V}ar\left( \sqrt{r}\left[ M_{N}\left( t\right) \right]
_{11}|Z\right) &=&\frac{2}{r}\sum_{i=1}^{r}\left( \frac{t}{t^{2}-s_{i}^{2}}%
\right) ^{2}-\frac{2}{r^{2}}\sum_{i\neq j}^{r}\left( \frac{t}{t^{2}-s_{i}^{2}%
}\right) \left( \frac{t}{t^{2}-s_{j}^{2}}\right) \\
&+&O\left( \frac{1}{r}\right)
\\
&\rightarrow &-2t^{2}\left( G^{\prime }\left( t^{2}\right) +\left[ G\left(
t^{2}\right) \right] ^{2}\right) .
\end{eqnarray*}%
In addition, Lemma \ref{speed_convergence} shows that 
\begin{equation*}
\sqrt{r}\left\{ \left[ M_{N}\left( t\right) \right] _{11}-\frac{1}{r}%
\sum_{i=1}^{r}\frac{t}{t^{2}-s_{i}^{2}}\right\}
\end{equation*}%
converges in distribution to the same variable as 
\begin{equation*}
\sqrt{r}\left\{ \left[ M_{N}\left( t\right) \right] _{11}-\mathbb{E}\frac{1}{%
r}\sum_{i=1}^{r}\frac{t}{t^{2}-s_{i}^{2}}\right\} ,
\end{equation*}%
that is, to a zero-mean Gaussian variable with variance $-2t^{2}\left(
G^{\prime }\left( t^{2}\right) +\left[ G\left( t^{2}\right) \right]
^{2}\right) .$

Similar argument holds for $\left[ M_{N}\left( t\right) \right] _{22}$: 
\begin{equation*}
\sqrt{r}\left\{ \left[ M_{N}\left( t\right) \right] _{22}-t\widetilde{G}%
_{N}\left( t^{2}\right) \right\}
\end{equation*}%
converges to a zero-mean Gaussian variable with variance 
\begin{equation*}
-2\beta ^{-1}t^{2}\left( \widetilde{G}^{\prime }\left( t^{2}\right) +\left[ 
\widetilde{G}\left( t^{2}\right) \right] ^{2}\right) .
\end{equation*}

Next, $\left[ M_{N}\left( t\right) +\theta ^{-1}\right] _{12}=u^{\ast
}\left( t^{2}I_{p}-ZZ^{\ast }\right) ^{-1}Zv,$ where $u$ and $v$ are unit
vectors, which are independent and uniformly distributed on $S^{p}$ and $S^{r},$ 
respectively. In appropriate coordinates, 
\begin{equation*}
\left[ M_{N}\left( t\right) +\theta ^{-1}\right] _{12}=\sum%
\limits_{i=1}^{p}u_{i}\frac{s_{i}}{t^{2}-s_{i}^{2}}v_{i}
\end{equation*}%
The expectation of this term is zero and for the conditional variance we
have the following sum:%
\begin{eqnarray*}
\mathbb{V}ar\left( \sqrt{r}\left[ M_{N}\left( t\right) +\theta ^{-1}\right]
_{12}|Z\right) &=&\frac{1}{p}\sum_{i=1}^{p}\left( \frac{s_{i}}{%
t^{2}-s_{i}^{2}}\right) ^{2} \\
&\rightarrow &\int \frac{s^{2}}{\left( t^{2}-s^{2}\right) ^{2}}d\widetilde{%
\mu }\left( s^{2}\right) \\
&=&-\left[ t^{2}\widetilde{G}^{\prime }\left( t^{2}\right) +\widetilde{G}%
\left( t^{2}\right) \right] .
\end{eqnarray*}%
\hfill $\square $

\begin{corollary}
\label{lemma_fluctuations_determinant} If $t>x_{2}$,
then the random variable 
\begin{equation*}
\sqrt{r}\left( \left| M_{N}\left( t\right)\right| - \left| 
\begin{array}{cc}
tG_{N}\left( t^{2}\right) & -\theta ^{-1} \\ 
-\theta ^{-1} & t\widetilde{G}_{N}\left( t^{2}\right)%
\end{array}%
\right|
\right).
\end{equation*}%
converges in distribution to a Gaussian random variable with the variance 
\begin{equation*}
\kappa ^{2}=t^{2}\left[ \widetilde{G}\left( t^{2}\right) \right] ^{2}\kappa
_{1}^{2}+t^{2}\left[ G\left( t^{2}\right) \right] ^{2}\kappa
_{2}^{2}+4\theta ^{-2}\tau ^{2},
\end{equation*}%
where $\kappa _{1},$ $\kappa _{2},$ and $\tau $ are as in Lemma \ref%
{lemma_M_fluctuations}.
\end{corollary}

Next, by Lemma \ref{lemma_reduction_to_M}, the largest singular value of $%
\widehat{A}$ satisfies the equation 
\begin{equation*}
\left| M_{N}\left( \widehat{\sigma }_{1}\right)\right| =0.
\end{equation*}%
Let 
\begin{equation*}
f_{N}\left( t\right) \df t^{2}\widetilde{%
G}_{N}\left( t^{2}\right) G_{N}\left( t^{2}\right) -\theta ^{-2},
\end{equation*}%
so that by Corollary \ref{lemma_fluctuations_determinant} 
\begin{equation*}
\left| M_{N}\left( t\right)\right| =f_{N}\left( t\right) +\frac{\kappa }{\sqrt{r}}%
W+o\left( r^{-1/2}\right) ,
\end{equation*}%
where $W$ is a standard Gaussian random variable. If $\sigma _{1}$ denote
the largest root of $f_{N}\left( t\right) =0,$ then this implies that 
\begin{eqnarray*}
\sqrt{r}\left( \widehat{\sigma }_{1}-\sigma _{1}\right) &\rightarrow &\left[
f_{N}^{\left( -1\right) }\right] ^{\prime }\left( 0\right) \kappa W= \\
&&\frac{1}{f_{N}^{\prime }\left( \sigma _{1}\right) }\kappa W,
\end{eqnarray*}%
Note that our estimator is $\widehat{\theta }=D_{N}\left( \widehat{\sigma }%
_{1}\right) ,$ where 
\begin{equation*}
D_{N}\left( t\right) =\sqrt{\frac{1}{t^{2}\widetilde{G_{N}}\left(
t^{2}\right) G_{N}\left( t^{2}\right) }}=\left( f_{N}\left( t\right) +\theta
^{-2}\right) ^{-1/2}.
\end{equation*}%
Hence 
\begin{eqnarray*}
\sqrt{r}\left( \widehat{\theta }-\theta \right) &\rightarrow &D_{N}^{\prime
}\left( \sigma _{1}\right) \frac{1}{f_{N}^{\prime }\left( \sigma _{1}\right) 
}\kappa W \\
&=&-\frac{1}{2}\left( f_{N}\left( \sigma _{1}\right) +\theta ^{-2}\right)
^{-3/2}\kappa W \\
&=&-\frac{1}{2}\theta ^{3}\kappa W.
\end{eqnarray*}%
\hfill $\square $

\section{Conclusion}

\label{conclusion} This paper is about the reduced-rank regression in the
multivariate response linear model. We found that if the number of responses
and predictors is large relative to the number of observations, then the
singular values of the OLS estimator of the coefficient matrix do not
converge to zero even if the true coefficient matrix is zero. The same
observation is true for the matrix of fitted responses. Instead, the
empirical distributions of singular values are converging to some limit
distributions that depend on how numerous the predictors and responses are
relative to the number of observations.

In addition, under the null hypothesis $A=0$ we found that the scaled
largest singular value for these matrices are distributed in the limit
according to the Tracy-Widom distribution. This fact can be used to test
whether the true coefficient matrix is zero and to estimate the rank of the model. In numerical simulations, we found that
one of these rank-selection algorithms compares favorably with the algorithm
suggested by Bunea, She, and Wegkamp in \cite{bunea_she_wegcamp11}.

In the case of the low-rank $A\neq 0$, we showed that the singular values of 
$A$ are detectable if and only if they exceed a certain threshold. If they
do, then the estimated coefficient matrix has singular values outside of the
support of the limit empirical distribution.

Finally, we showed that consistent estimators of the true singular values
can be obtained by shrinking the outlier singular values of $\widehat{A}$ or $\widehat{Y}$
appropriately. 
We found that the estimation based on singular values of $\widehat{Y}$ is preferable to the estimation based on that 
of $\widehat{A}.$  We have also proved a CLT for the 
asymptotic distribution of one of the estimators. 

\appendix

\section{The limiting distribution of singular values for the coefficient matrix estimator}

\label{section_limiting_distribution}

In this section we will prove the second part of Theorem \ref{Theorem_main} by using the technique of $S$-transforms from free probability. (This technique was introduced in \cite{voiculescu87} and generalized in \cite{bercovici_voiculescu93} in order to study the spectra of products of infinite-dimensional operators.) Let the moments of
a square $N$-by-$N$ random matrix $A$ be defined as $m_{k}\df N^{-1}\mathrm{Tr}\left( A^{k}\right),$ and let the generating function for $m_k$ be denoted as 
\begin{equation*}
M_{A}\left( z\right) \df\sum_{k=1}^{\infty }m_{k}z^{k}.
\end{equation*}%
The Stieltjes transform of the empirical eigenvalue distribution is 
\begin{equation*}
G_{A}\left( z\right) = \sum_{k=0}^{\infty }\frac{m_{k}}{z^{k+1}}=\frac{1}{z}+\frac{1}{z}M_{A}\left( \frac{1}{z}%
\right) .
\end{equation*}%
and Voiculescu's $S$-transform is 
\begin{equation*}
S_{A}\left( u\right) \df\frac{u+1}{u}M_{A}^{\left( -1\right) }\left( u\right)
,
\end{equation*}%
where $M_{A}^{\left( -1\right) }$ is the functional inverse of $M_A$.

If $X_{N}$ is an $N$-by-$p$ random matrix with the standard Gaussian
entries, and $N/p\to 1+\lambda $, then it is known that as $%
N\to \infty ,$ the moments of matrices $\left( X_{N}^{\ast
}X_{N}\right) /p$ converge to the moments of the Marchenko-Pastur
distribution with parameter $1+\lambda .$ This distribution has the $S$-transform $S_{MP}\left( u\right) =1/\left(
u+1+\lambda \right) .$

\begin{lemma}
\label{lemma_inv_MP}Let $\xi $ be a random variable with the
Marchenko-Pastur distribution with parameter $1+\lambda \geq 1$. Then, the $%
S $-transform of $\xi ^{-1}$ is 
\begin{equation*}
S_{\xi ^{-1}}\left( u\right) =\lambda -u.
\end{equation*}
\end{lemma}

\textbf{Proof:} Let $F\left( \alpha \right) \df\alpha -\sqrt{\alpha ^{2}-1}$.
Then, a calculation of integrals gives the moments of $\xi ^{-1}$: 
\begin{equation*}
m_{k}=\int t^{-k}dF_{MP}\left( t\right) =\sqrt{1+\lambda }\left( \frac{-1}{2%
\sqrt{1+\lambda }}\right) ^{k}\frac{1}{k!}F^{\left( k\right) }\left( \frac{%
2+\lambda }{2\sqrt{1+\lambda }}\right) ,
\end{equation*}%
and therefore, 
\begin{equation*}
M_{\xi ^{-1}}\left( z\right) =\frac{1}{2}\left( \lambda -z-\sqrt{\left( \lambda
+2-z\right) ^{2}-4\left( 1+\lambda \right) }\right) .
\end{equation*}%
Hence, 
\begin{equation*}
M_{\xi ^{-1}}^{\left( -1\right) }\left( u\right) =\frac{u}{u+1}\left( \lambda -u\right) ,
\end{equation*}%
and 
\begin{equation*}
S_{\xi ^{-1}}\left( u\right) =\lambda -u.
\end{equation*}%
\hfill $\square $

\begin{lemma}
\label{lemma_shift}Let $A$ be an $N$-by-$p$ matrix, $B$ be an $p$-by-$p$
matrix, $Z=ABA^{\ast },$ and $\widetilde{Z}=BA^{\ast }A.$ Then, the
Stieltjes and $S$-transforms of matrices $Z$ and $\widetilde{Z}$ are related
as follows:%
\begin{eqnarray*}
G_{Z}\left( t\right) &=&\left( 1-\frac{p}{N}\right) \frac{1}{t}+\frac{p}{N}%
G_{\widetilde{Z}}\left( t\right) \\
S_{Z}\left( u\right) &=&\frac{u+1}{u+\frac{p}{N}}S_{\widetilde{Z}}\left( 
\frac{N}{p}u\right) .
\end{eqnarray*}
\end{lemma}

\textbf{Proof: }The moments of matrices $Z$ and $\widetilde{Z}$ are related
as follows: 
\begin{eqnarray*}
m_{k} &=&\frac{1}{N}Tr\left( \left[ ABA^{\ast }\right] ^{k}\right) \\
&=&\frac{p}{N}\frac{1}{m}Tr\left( \left[ BA^{\ast }A\right] ^{k}\right) =%
\frac{p}{N}\widetilde{m}_{k}.
\end{eqnarray*}%
The rest follows from the definitions of the Stieltjes and S- transforms. 
 \hfill $ \square $

\begin{lemma}
Let $A_{N}\df\left( X_{N}^{\ast }X_{N}\right) ^{-1}X_{N}^{\ast }U_{N}.$  
As $N\to \infty ,$ the $S$-transform of $\frac{p}{r}A_{N}^{\ast }A_{N}$
converges to 
\begin{equation*}
\frac{\lambda - u/\beta }{u+\beta}.
\end{equation*}
\end{lemma}

\textbf{Proof:} By applying Lemma \ref{lemma_shift} to matrices \\ 
$\widetilde{Z}=p\left( X_{N}^{\ast }X_{N}\right) ^{-1}=p\left( X_{N}^{\ast }X_{N}\right)
^{-2}\left( X_{N}^{\ast }X_{N}\right) $ and $Z=pX_{N}\left( X_{N}^{\ast
}X_{N}\right) ^{-2}X_{N}^{\ast },$ we find that 
\begin{equation*}
S_{Z}\left( u\right) =\frac{u+1}{u+\frac{p}{N}}S_{\widetilde{Z}}\left( \frac{%
N}{p}u\right) .
\end{equation*}%
If $\lambda >1$ then the moments of $p\left( X_{N}^{\ast }X_{N}\right) ^{-1}$
converge to the moments of the inverse Marchenko-Pastur distribution and
therefore (by Lemma \ref{lemma_inv_MP}) $S_{\widetilde{Z}}\left( \frac{N}{p}%
u\right) $ converges to $\lambda -(1+\lambda )u.$ Therefore, as $%
N\rightarrow \infty ,$ 
\begin{equation*}
S_{Z}\left( u\right) \rightarrow \phi_1(u)\equiv\frac{u+1}{u+\frac{1}{1+\lambda }}\left(
\lambda -\left( 1+\lambda \right) u\right) .
\end{equation*}

Next, the matrices $r^{-1}\left( U_{N}^{\ast }U_{N}\right) $ converge in
distribution to a multiple of the Marchenko-Pastur distribution with
parameter $\mu ,$ and one compute (for example, by Lemma \ref{lemma_shift})
that the $S$ transform of $r^{-1}\left( U_{N}U_{N}^{\ast }\right) $
converges to $\phi_2(u)\equiv\left( \mu u+1\right) ^{-1}.$

Next, we use the facts that $X_{N}$ and $Y_{N}$ are asymptotically free, and
that the $S$-transform of the product of free variables is the product of
their $S$-transforms. Hence the $S$-transform of the matrix 
\begin{equation*}
\frac{p}{r}X_{N}\left( X_{N}^{\ast }X_{N}\right) ^{-2}X_{N}^{\ast
}U_{N}U_{N}^{\ast }
\end{equation*}%
converges to $\phi_1(u) \phi_2(u).$
Another application of Lemma \ref{lemma_shift} shows that the $S$-transform of\\ 
$\frac{p}{r}U_{N}^{\ast }X_{N}\left( X_{N}^{\ast }X_{N}\right) ^{-2}X_{N}^{\ast }U_{N}$ converges to 
\begin{equation*}
\left(\lambda -\left( \frac{1+\lambda }{\mu 
}\right) u \right)\left(u+\frac{\mu }{1+\lambda }\right)^{-1}.
\end{equation*}%
\hfill $\square $

Now, by inverting$\left[ u/\left( u+1\right) \right] S\left( u\right) ,$ we
calculate $M\left( z\right) $ and then the Stieltjes transform $G(z)$ for the limit of the
variables $\frac{p}{r}A_{N}^{\ast }A_{N}:$%
\begin{equation}
G\left( z\right) =\frac{\beta }{2}\frac{1-\beta +\left( \lambda +\frac{2}{%
\beta }\right) z-\sqrt{\left[ \left( 1+\beta \right) -\lambda z\right]
^{2}-4\left( z+\beta \right) }}{z\left( z+\beta \right) },
\label{Stieltjes_P_lambda_beta}
\end{equation}%
where $\beta \df\mu /(1+\lambda).$
From this we can extract the density function and the support of the
limiting probability measure for $\frac{p}{r}A_{N}^{\ast
}A_{N}. $ The limit distribution for $A_{N}^{\ast}A_{N}$ is obtained by scaling this distribution by $%
\beta ^{-1}.$ \hfill $\square $

\bibliographystyle{plain}
\bibliography{comtest}

\end{document}